\documentclass[12pt]{amsart}

\usepackage[utf8]{inputenc}
\usepackage[T1,T2A]{fontenc}
\usepackage[english]{babel}

\usepackage{amssymb}
\usepackage{amsmath}
\usepackage{amsthm}
\usepackage{enumitem}
\usepackage{array}
\usepackage[hidelinks]{hyperref}
\usepackage{tikz}
\usetikzlibrary{calc}
\usepackage{caption}
\usepackage{subcaption}
\usepackage{float}

\newcommand{\CC}{\mathbb{C}}
\newcommand{\ZZ}{\mathbb{Z}}
\newcommand{\QQ}{\mathbb{Q}}
\newcommand{\PP}{\mathbb{P}}
\newcommand{\AAA}{\mathbb{A}}

\newcommand{\FF}{\mathfrak{F}}
\renewcommand{\phi}{\varphi}
\newcommand{\eps}{\varepsilon}

\newcommand{\rk}{\operatorname{rk}}
\newcommand{\Hom}{\operatorname{Hom}}
\newcommand{\cone}{\operatorname{cone}}
\newcommand{\Bl}{\operatorname{Bl}}
\newcommand{\Star}{\operatorname{Star}}
\newcommand{\Supp}{\operatorname{Supp}}
\newcommand{\Ima}{\operatorname{Im}}
\newcommand{\GL}{\operatorname{GL}}
\newcommand{\SL}{\operatorname{SL}}
\newcommand{\diag}{\operatorname{diag}}
\newcommand{\pf}{\operatorname{pf}}

\newcommand{\pairing}[2]{\langle#1,\,#2\rangle}

\newtheorem{theorem}{Theorem}[section]
\newtheorem{proposition}[theorem]{Proposition}
\newtheorem{lemma}[theorem]{Lemma}
\newtheorem{corollary}[theorem]{Corollary}

\theoremstyle{remark}
\newtheorem{remark}[theorem]{Remark}

\theoremstyle{definition}
\newtheorem{definition}[theorem]{Definition}
\newtheorem{example}{Example}[section]

\begin{document}

\date{}
\title[Combinatorics of equivariant blowups of spherical varieties]{Combinatorics of equivariant blowups of spherical varieties}

\author{Daniil Shunin}
\address{
Lomonosov Moscow State University, Moscow, Russia
\newline\indent HSE University, Moscow, Russia}
\email{sshunindaniil@gmail.com}
\thanks{This article is an output of a research project implemented as part of the Basic Research Program at HSE University}


\begin{abstract}
We explicitly describe the colored fan of the blowup of a smooth spherical variety $X$ along a smooth irreducible $G$-stable subvariety $Y$. The answer is expressed in terms of the colored cone of $Y$ and the localization of the valuation cone at $Y$, and generalizes the star subdivision of fans of toric varieties. We also consider several examples illustrating the main result.
\end{abstract}

\maketitle


\section{Introduction}

Let $G$ be a connected reductive complex algebraic group. A normal irreducible $G$-variety is called \textit{spherical} if a Borel subgroup $B$ of $G$ acts on it with an open orbit. An embedding of a spherical homogeneous space $G/H$ into a normal $G$-variety $X$ as an open dense $G$-orbit is called a \textit{spherical embedding}. The spherical embeddings of a given homogeneous space $G/H$ are classified by the Luna--Vust theory~\cite{Kn_LV}.

Several birational invariants are attached to the space $G/H$, such as the lattice $\mathcal{M}=\mathcal{M}_{G/H}$ of weights of rational $B$-eigenfunctions and the set of \textit{colors} $\mathcal{D}=\mathcal{D}_{G/H}$, that is, the set of prime $B$-stable divisors on $G/H$. To each spherical embedding $G/H\hookrightarrow X$ one then assigns its \textit{colored fan} $\FF(X)$, which consists of finitely many \textit{colored cones}, that is, pairs $(\mathcal{C},\mathcal{F})$, where $\mathcal{F}\subset\mathcal{D}$ and $\mathcal{C}$ is a pointed polyhedral cone in the rational vector space $\mathcal{N}_\QQ$ spanned by the lattice $\mathcal{N}=\mathcal{N}_{G/H}=\Hom_\ZZ(\mathcal{M},\ZZ)$. By the Luna--Vust theory, the map $G/H\subset X\mapsto\FF(X)$ is a bijection.

Many properties of a spherical variety and of the geometric objects associated with it can be read off from the corresponding colored fan (adjacency of orbits, local factoriality, smoothness, the Picard group, the canonical class, etc.). For example, the $G$-orbits in $X$ (there are finitely many of them, see~\cite{Akh}) are in bijection with the colored cones of $\FF(X)$. Since an orbit $O$ and its closure $Y=\overline{O}\subset X$ determine each other, it will be more convenient for us to speak of the bijection $Y\mapsto(\mathcal{C}(Y),\mathcal{F}(Y))$ between the irreducible closed $G$-stable subvarieties of $X$ and the colored cones of $\FF(X)$.

Let $Y\subset X$ be a closed irreducible $G$-stable subvariety of a spherical variety $X\supset G/H$. The blowup $\Bl_YX$ of $X$ along $Y$ carries an action of $G$ and, being birational to $X$, contains $G/H$. The variety $\Bl_YX$ may be non-normal, but its normalization $X'$ is a spherical embedding of $G/H$. This leads to the problem of describing the colored fan $\FF(X')$, given the fan $\FF(X)$ of the original variety, the colored cone $(\mathcal{C}(Y),\mathcal{F}(Y))$ corresponding to $Y$, and the additional combinatorial data of the open $G$-orbit that are standard in the theory of spherical varieties.

For instance, if $X$ is a toric variety, that is, a spherical variety for an algebraic torus $T$, the solution for the blowup of a smooth $T$-fixed point $x\in X$ is well known. Namely, we may assume that $T$ acts effectively on the dense orbit; then $X$ corresponds to an ordinary fan in the space $N_\QQ=N\otimes_\ZZ\QQ$, where $N$ is the lattice of one-parameter subgroups of $T$. The toric subvariety having $x$ as its only closed orbit corresponds to a rational cone $\cone(u_1,\dots,u_k)$ with $u_i\in N$, and if $X$ is smooth at $x$, then $k=\rk N$ and $u_1,\dots,u_k$ form a basis of $N$. The fan of $\Bl_xX$ is obtained from the fan of $X$ by the so-called \textit{star subdivision} at the vector $u_0=u_1+\dots+u_k$. The blowup of a smooth stable irreducible subvariety of a toric variety is described similarly (see, e.g., \cite[\S\,3.3, ``Refinements of Fans and Blowups'']{CLS}).

The most general answer to this problem was given by M.~Brion~\cite{Br}: he described the colored fan of the (normalized) blowup $X'$ of an invariant ideal sheaf $\mathcal{I}$ (that is, a coherent sheaf of submodules of $\mathcal{O}_X$) or, more generally, of a \textit{fractional ideal sheaf} (that is, a coherent sheaf of $\mathcal{O}_X$-submodules of the sheaf of rational functions on $X$) on a spherical variety $X$. Namely, from the sheaf $\mathcal{I}$ one builds a special piecewise linear function $\operatorname{ord}\mathcal{I}$ on the support of the colored fan $\FF(X)$ in $\mathcal{N}_\QQ$ (see \S\,\ref{s: equiv_morphisms}), and the domains of linearity of this function determine how $\FF(X)$ must be subdivided to obtain~$\FF(X')$.

The one-dimensional colored cones without colors in the colored fan of a spherical variety correspond to its $G$-stable prime divisors. Those one-dimensional cones of $\FF(X')$ that do not occur in $\FF(X)$ correspond to the irreducible components of the exceptional divisor $E=\pi^{-1}(Y)$ of the blowup $\pi\colon X'\to X$. For smooth $X$ and $Y$, the exceptional divisor $E\subset\Bl_YX$ is prime and isomorphic to the total space of the projectivization $\PP(N_{Y/X})$ of the normal bundle of $Y$ in $X$. In general, $E$ may have several components, but if $X$ is smooth, only one of them maps dominantly onto $Y$. In the toric case with $X$ smooth and $Y=\mathrm{pt}$ considered above, for instance, the exceptional divisor corresponds to the ray spanned by~$u_0$.

So, in the case of smooth $X$ and $Y$ a single ray gets glued into the fan $\mathfrak{F}(X)$. This ray enters Brion's description indirectly as the ray that separates the domains of linearity of $\operatorname{ord}\mathcal{I}$, where $\mathcal{I}$ is the ideal sheaf of $Y$. To determine the colors of $\mathfrak{F}(X')$ one also has to find the linear forms which constitute the piecewise linear function $\operatorname{ord}\mathcal{I}$. We suggest a different description of $\mathfrak{F}(X')$ in terms of spherical roots which does not involve the function $\operatorname{ord}\mathcal{I}$. Instead, we make use of smoothness by applying Luna's \'etale slice theorem and reducing the problem to spherical modules.

The goal of this paper is to describe explicitly the colored fan $\FF(X')$ of the blowup $X'$ of a smooth spherical variety $X$ along a smooth irreducible $G$-stable subvariety $Y$. The paper is organized as follows. In \S\,\ref{s: preliminaries} we collect the facts needed to prove the main theorem: the combinatorial description of spherical embeddings and invariant valuations on spherical varieties (\S\,\ref{s: combinatorics}), the description of equivariant morphisms in terms of colored fans (\S\,\ref{s: equiv_morphisms}), the local structure theorem, with an emphasis on the behavior of birational invariants under localization (\S\,\ref{s: local_structure}), and a spherical analogue of the star subdivision from the theory of toric varieties (\S\,\ref{s: star}). In \S\,\ref{s: main_result} we state and prove the main theorem (Theorem~\ref{thm: main}), which gives an explicit description of $\FF(X')$. Section~\ref{s: examples} contains examples illustrating the main result.

The author is grateful to D.\,A.~Timashev and V.\,S.~Zhgoon for their atention to the paper.

\section{Preliminaries}
\label{s: preliminaries}

Let $G$ be a connected complex reductive algebraic group. We fix a Borel subgroup $B\subset G$ and a maximal torus $T\subset B$. As usual, we identify the character groups of $B$ and $T$. Denote by $R$ the root system of $G$ with respect to $T$, and by $S\subset R$ the set of simple roots determined by $B$. Each subset $S_*\subset S$ corresponds to a unique standard parabolic subgroup $P_{S_*}$, namely, the minimal parabolic subgroup $P\supset B$ containing the root subgroups $U_{-\alpha}$ for all $\alpha\in S_*$. We put $P_{\{\alpha\}}=:P_\alpha$.

\subsection*{Notation and conventions}

\begin{itemize}[leftmargin=*,itemsep=0pt]
\item[] $\mathcal{C}^\circ$ is the relative interior of a polyhedral cone $\mathcal{C}$.
\item[] $\overline{Y}$ is the Zariski closure of a subset $Y$ in a variety $X$.
\item[] $\CC[X]$ is the algebra of regular functions on a variety $X$.
\item[] $\CC(X)$ is the field of rational functions on an irreducible variety $X$.
\item[] $\cone(A)$ is the polyhedral cone generated by a finite set $A$.
\item[] $A^{(G)}$ is the set of $G$-eigenvectors in a $G$-module $A$.
\item[] $A^{(G)}_\lambda$ is the set of $G$-eigenvectors of weight $\lambda$ in a $G$-module $A$.
\end{itemize}
The base field is that of complex numbers $\CC$. All varieties and groups are assumed to be algebraic, and all topological terms refer to the Zariski topology.

\subsection{Combinatorics of spherical varieties}
\label{s: combinatorics}

Let $G/H$ be a spherical homogeneous space. Denote by $\mathcal{M}=\mathcal{M}_{G/H}$ the lattice of weights of rational $B$-eigenfunctions, and by $\mathcal{N}=\mathcal{N}_{G/H}=\Hom(\mathcal{M},\ZZ)$ the dual lattice. Each weight $\lambda\in\mathcal{M}$ corresponds to a $B$-eigenfunction $f_\lambda$, which is unique up to a scalar factor. A $\QQ$-valued discrete valuation $v$ of the field $\CC(G/H)$ defines a linear function $\overline{v}\in\mathcal{N}_\QQ:=\mathcal{N}\otimes_\ZZ\QQ$ on $\mathcal{M}$ by $\overline{v}(\lambda)=v(f_\lambda)$. Restricted to the set $\mathcal{V}$ of $G$-invariant discrete valuations, this assignment is an embedding $\mathcal{V}\hookrightarrow\mathcal{N}_\QQ$~\cite{Kn_LV}. It is known that the image of $\mathcal{V}$ is a cosimplicial polyhedral cone, called the \textit{valuation cone}, which we also denote by $\mathcal{V}$. When no confusion can arise, we use the same symbol $v$ both for a $G$-invariant valuation and for its image in $\mathcal{N}_\QQ$. The primitive outer normal vectors of the cone $\mathcal{V}$ are linearly independent elements of $\mathcal{M}$; they are called \textit{spherical roots}, and we denote the set of them by $\Sigma=\Sigma_{G/H}\subset\mathcal{M}$. The linear part $\mathcal{V}\cap(-\mathcal{V})$ of the valuation cone is denoted by $l\mathcal{V}$.

The $B$-stable prime divisors on $G/H$ are called \textit{colors}; we denote the set of them by $\mathcal{D}=\mathcal{D}_{G/H}$. To each color $D\in\mathcal{D}$ we assign the set $\varsigma(D)\subset S$ of those simple roots $\alpha$ for which $P_\alpha\cdot D\neq D$. We say that such roots $\alpha$ \textit{move} the color $D$. For $\alpha\in S$, denote the set of colors moved by $\alpha$ by
\begin{equation*}
\mathcal{D}(\alpha)=\{D\in\mathcal{D}\mid\alpha\in\varsigma(D)\}.
\end{equation*}
We also define the map $\rho\colon\mathcal{D}\to\mathcal{N}$, where $\rho(D)$ is the linear function whose value at a weight $\lambda\in\mathcal{M}$ is the order of $f_\lambda$ along the color $D$.

Spherical embeddings are classified by colored fans.

\begin{definition}
A \textit{colored cone} is a pair $(\mathcal{C},\mathcal{F})$, where $\mathcal{F}\subset\mathcal{D}$, the set $\rho(\mathcal{F})$ does not contain $0$, and $\mathcal{C}\subset\mathcal{N}_\QQ$ is a pointed cone generated by $\rho(\mathcal{F})$ and finitely many elements of $\mathcal{V}$, whose relative interior has a nonempty intersection with $\mathcal{V}$.

A \textit{colored face} of a colored cone $(\mathcal{C},\mathcal{F})$ is a colored cone $(\mathcal{C}',\mathcal{F}')$ such that $\mathcal{C}'$ is a face of $\mathcal{C}$, and $\mathcal{F}'\subset\mathcal{F}$ consists of the colors whose images under the map $\rho$ lie in $\mathcal{C}'$.
\end{definition}

\begin{definition}
A \textit{colored fan} is a finite collection $\FF$ of colored cones which is closed under taking colored faces and such that each $v\in\mathcal{V}$ belongs to the relative interior of at most one member of the collection.
\end{definition}

Let $X\supset G/H$ be a spherical embedding. To each closed irreducible $G$-stable subvariety $Y\subset X$ we assign the set $\mathcal{F}(Y)=\{D\in\mathcal{D}\mid Y\subset\overline{D}\}$ and the cone $\mathcal{C}(Y)\subset\mathcal{N}_\QQ$ generated by $\rho(\mathcal{F}(Y))$ and by all invariant valuations $v_D\in\mathcal{V}$ corresponding to the $G$-stable prime divisors $D$ of $X$ that contain $Y$. The pair $(\mathcal{C}(Y),\mathcal{F}(Y))$ is a colored cone, and these colored cones, taken over all closed irreducible $G$-stable subvarieties $Y\subset X$, form the colored fan $\FF(X)$~\cite{Kn_LV}. The maximal cones $\mathcal{C}(Y)$ of $\FF(X)$ correspond to the open $G$-stable \textit{simple} spherical subvarieties of $X$ (those with a unique closed orbit $Y$), which cover $X$.

We also recall some facts about centers of invariant valuations on spherical varieties. Let $X$ be an arbitrary irreducible normal variety, and $v$ a discrete $\QQ$-valued valuation of the field $\CC(X)$. The valuation ring $R_v=\{f\in\CC(X)\mid v(f)\geq0\}$ is local, with maximal ideal $\mathfrak{m}_v=\{f\in\CC(X)\mid v(f)>0\}$.

\begin{definition}
An irreducible subvariety $Z\subset X$ is called the \textit{center} of the valuation $v$ if $\mathcal{O}_{X,Z}\subset R_v$ and $\mathfrak{m}_{X,Z}\subset\mathfrak{m}_v$, where $\mathcal{O}_{X,Z}$ is the local ring at the generic point of $Z$ and $\mathfrak{m}_{X,Z}$ is its maximal ideal.
\end{definition}

A valuation has at most one center. Clearly, if a group acts on $X$ and the valuation is invariant, then its center, if it exists, is also stable.

\begin{proposition}[{\cite[Theorem 2.5]{Kn_LV}}]
\label{prop: centers}
Let $X$ be a spherical variety, and let $v\in\mathcal{V}_X$.
\begin{enumerate}[label=(\arabic*),itemsep=0pt]
\item The valuation $v$ has a center if and only if $v\in\mathcal{C}$ for some colored cone $(\mathcal{C},\mathcal{F})$ of the fan $\FF(X)$.
\item A closed irreducible $G$-stable subvariety $Y$ is the center of $v$ if and only if $v$ lies in $\mathcal{C}(Y)^\circ$.
\end{enumerate}
\end{proposition}

\begin{remark}
\label{rm: inclusion_reversing}
Note that the correspondence between colored cones and closed $G$-stable subvarieties reverses inclusions: $Y_1\subset Y_2$ if and only if $(\mathcal{C}(Y_2),\mathcal{F}(Y_2))$ is a colored face of $(\mathcal{C}(Y_1),\mathcal{F}(Y_1))$.
\end{remark}

Finally, we recall a characterization of the valuation cone in the affine case. Let $X$ be an affine spherical variety. Denote by $\CC[X]_{(\lambda)}$ the isotypic component of the $G$-module $\CC[X]$ corresponding to the highest weight $\lambda$ (that is, the sum of all simple $G$-submodules of highest weight $\lambda$). The variety $X$ is spherical if and only if the $G$-module $\CC[X]$ is \textit{multiplicity free}, that is, in the decomposition
\begin{equation*}
\CC[X]=\bigoplus_{\lambda\in\mathcal{M}_X^+}\CC[X]_{(\lambda)}
\end{equation*}
into isotypic components, every $G$-module $\CC[X]_{(\lambda)}$ is simple. Here $\mathcal{M}_X^+$ is the semigroup of weights of those functions $f_\lambda$ that are regular on $X$. For $\lambda,\mu\in\mathcal{M}_X^+$, the product of the corresponding components decomposes into simple submodules as
\begin{equation*}
\CC[X]_{(\lambda)}\cdot\CC[X]_{(\mu)}=\CC[X]_{(\lambda+\mu)}\oplus\bigoplus_{i=1}^k\CC[X]_{(\lambda+\mu-\beta_i)},
\end{equation*}
where the $\beta_i$ are some sums of simple roots of $G$ (called \textit{tails}). Put
\begin{equation*}
\mathcal{T}_X=\cone(\beta_i\mid\lambda,\mu\in\mathcal{M}_X^+).
\end{equation*}

\begin{proposition}[{\cite[\S\,21.1]{Ti}}]
\label{prop: valuation_cone_affine}
Then
\begin{equation*}
\mathcal{V}_X=-\mathcal{T}_X^\vee,
\end{equation*}
where $\mathcal{T}_X^\vee=\{v\in\mathcal{N}_\QQ\mid\pairing{\mathcal{T}_X}{v}\geq0\}$ denotes the dual cone.
\end{proposition}

\subsection{Equivariant morphisms}
\label{s: equiv_morphisms}

Let us recall how equivariant morphisms between spherical varieties are described in terms of their colored fans.

An equivariant morphism $\phi\colon G/H\to G/H'$ of spherical homogeneous spaces induces homomorphisms of the corresponding lattices
\begin{equation*}
\phi^*\colon\mathcal{M}_{G/H'}\to\mathcal{M}_{G/H},\qquad\phi_*\colon\mathcal{N}_{G/H}\to\mathcal{N}_{G/H'},
\end{equation*}
and it is known that $\phi_*(\mathcal{V}_{G/H})=\mathcal{V}_{G/H'}$. Denote by $\mathcal{D}_\phi$ the set of all colors $D\in\mathcal{D}_{G/H}$ whose images are dense in $G/H'$. We also have a surjective map
\begin{equation*}
\phi_*\colon\mathcal{D}_{G/H}\backslash\mathcal{D}_\phi\to\mathcal{D}_{G/H'}.
\end{equation*}

\begin{definition}
Let $(\mathcal{C},\mathcal{F})$ and $(\mathcal{C}',\mathcal{F}')$ be the colored cones of some simple spherical embeddings of the homogeneous spaces $G/H$ and $G/H'$, respectively. We say that $(\mathcal{C},\mathcal{F})$ \textit{maps} to $(\mathcal{C}',\mathcal{F}')$ if
\begin{enumerate}[label=(\arabic*),itemsep=0pt]
\item $\phi_*(\mathcal{C})\subset\mathcal{C}'$,
\item $\phi_*(\mathcal{F}\backslash\mathcal{D}_\phi)\subset\mathcal{F}'$.
\end{enumerate}
\end{definition}

\begin{definition}
Let $\FF$ and $\FF'$ be the colored fans of some spherical embeddings of the homogeneous spaces $G/H$ and $G/H'$, respectively. We say that $\FF$ \textit{maps} to $\FF'$ if each colored cone of the fan $\FF$ maps to some colored cone of the fan $\FF'$.
\end{definition}

\begin{definition}
The \textit{support} $\Supp\FF$ of a colored fan $\FF$ is the following subset of the space~$\mathcal{N}_\QQ$:
\begin{equation*}
\Supp\FF=\mathcal{V}\cap\Bigl(\bigcup_{(\mathcal{C},\mathcal{F})\in\FF}\mathcal{C}\Bigr).
\end{equation*}
\end{definition}

\begin{theorem}
\label{thm: morphisms}
Let $\phi\colon G/H\to G/H'$ be an equivariant morphism between spherical homogeneous spaces, and let $G/H\subset X$ and $G/H'\subset X'$ be two spherical embeddings. Then $\phi$ extends to a morphism $X\to X'$ if and only if the colored fan $\FF(X)$ maps to the colored fan $\FF(X')$. In this case, the extension of $\phi$ is a proper morphism if and only if $\Supp\FF(X)=\phi_*^{-1}(\Supp\FF(X'))\cap\mathcal{V}_X$.
\end{theorem}

We will only need a very special case of this description.

\begin{corollary}
\label{cor: proper}
If $X$ and $X'$ are two spherical embeddings of the same homogeneous space $G/H$, then a dominant $G$-equivariant morphism $\phi\colon X\to X'$ is proper if and only if the supports of the fans $\FF(X)$ and $\FF(X')$ coincide.
\end{corollary}

\subsection{Local structure theorem}
\label{s: local_structure}

Let us recall the local structure theorem for spherical varieties. Let $Y\subset X$ be a closed irreducible $G$-stable subvariety, and let $(\mathcal{C},\mathcal{F})$ be the corresponding colored cone. Denote $X_Y=\{x\in X\mid\overline{Gx}\supset Y\}$ and put
\begin{equation*}
X^\circ=X_Y\backslash\bigcup_{D\in\mathcal{D}\backslash\mathcal{F}}\overline{D}.
\end{equation*}
Let $P$ be the stabilizer of the set $X^\circ$, and let $P=P_u\leftthreetimes L$ be its Levi decomposition.

\begin{theorem}[{\cite[Theorem 15.17]{Ti}}]
\label{thm: local_structure}
There exists an $L$-spherical affine subvariety $Z\subset X^\circ$ such that:
\begin{enumerate}[label=(\arabic*),itemsep=0pt]
\item There is a $P$-equivariant isomorphism
\begin{equation*}
X^\circ\cong P*_LZ\cong P_u\times Z.
\end{equation*}
\item Denote by $L_0$ the kernel of the action of $L$ on $Y\cap Z$. Then $L_0\supset[L,L]$ and $Y\cap Z\cong L/L_0$.
\item There exists an $L_0$-stable spherical affine subvariety $Z_0\subset Z$ that intersects $Y\cap Z$ transversally at its unique fixed point, and $Z\cong L*_{L_0}Z_0$.
\end{enumerate}
Moreover, $\mathcal{M}_Z=\mathcal{M}_X$ and $\rk\mathcal{M}_{Z_0}=\dim\mathcal{C}$.
\end{theorem}

We need to pay special attention to birational invariants and colored fans of the varieties $X$, $Z$ and $Z_0$. Let $(\mathcal{C}_Z,\mathcal{F}_Z)$ and $(\mathcal{C}_{Z_0},\mathcal{F}_{Z_0})$ be the colored cones corresponding to the unique closed orbits of $Z$ and of $Z_0$, respectively. First, the map $f\mapsto f|_Z$ is a bijection between the $B$-eigenfunctions on $X$ and the $(B\cap L)$-eigenfunctions on $Z$, so $\mathcal{M}_Z=\mathcal{M}_X$. Next, there is a bijection between the $B$-stable prime divisors on $X$ containing $Y$ and the $(B\cap L)$-stable prime divisors on $Z$ containing the closed orbit, and corresponding divisors have the same image in the lattice $\mathcal{N}_X=\mathcal{N}_Z$. Hence $\mathcal{C}_Z=\mathcal{C}$.

In~\cite{Ga}, a localization procedure is described that allows one to compute the birational invariants of the $L$-variety $Z$ combinatorially. We need only the description of the valuation cone $\mathcal{V}_Z$, namely, $\mathcal{V}_Z=-\Sigma_Z^\vee$, where
\begin{equation*}
\Sigma_Z=\Sigma_X\cap\langle S_\mathcal{F}\rangle_\QQ,\qquad S_\mathcal{F}=\{\alpha\in S\mid\mathcal{D}(\alpha)\subset\mathcal{F}\}.
\end{equation*}
Note that the set of simple roots $S_\mathcal{F}\subset S$ is such that $P=P_{S_\mathcal{F}}$.

Let $M(T)$ denote the weight group of any torus $T$, and set $N(T)=\Hom(M(T), \ZZ)$. Also set $L_0^\circ$ to be the identity component of $L_0$.

\begin{lemma}
\label{lemma: localization}
There is a short exact sequence of lattices
$$
0\to M(L/L_0)\to \mathcal{M}_Z\to\mathcal{M}_{Z_0}\to 0
$$
induced by the natural projection $Z=L*_{L_0}Z_0\to L/L_0$ and the restriction of weights from $B\cap L$ to $B\cap(L_0)^\circ$.
The dual embedding $\mathcal{N}_{Z_0}\hookrightarrow\mathcal{N}_Z$ is saturated, i.e. the abelian group $\mathcal{N}_Z/\mathcal{N}_{Z_0}$ is torsion-free. On the level of rational spaces and cones we have
\begin{align*}
(\mathcal{N}_Z)_\QQ &= N(L/L_0)_\QQ\times (\mathcal{N}_{Z_0})_\QQ, \\
\mathcal{V}_Z &= N(L/L_0)_\QQ\times\mathcal{V}_{Z_0}, \\
l\mathcal{V}_Z &\cap (\mathcal{N}_{Z_0})_\QQ = l\mathcal{V}_{Z_0}
\end{align*}
(under some splitting of the dual short exact sequence over $\QQ$).
In particular, $\mathcal{V}_Z=\mathcal{V}_{Z_0}+l\mathcal{V}_Z$.
\end{lemma}
\begin{proof}
Since $[L, L]\subset L_0^\circ$ we may choose a torus $T'$ in the connected center of $L$ such that $L=T' L_0^\circ$ and $K=T'\cap L_0$ is finite. Then $Z\cong T' *_K Z_0$.

Let $\lambda\in\mathcal{M}_{Z_0}$. As $Z_0$ is $L_0^\circ$-spherical and $K$ centralizes $B_0=B\cap L_0^\circ$, $K$ preserves the one-dimensional space of rational $B_0$-eigenfunctions of weight $\lambda$. Thus, there is a character $\chi\in M(K)$ such that
$$
f_\lambda(k z) = \chi(k) f_\lambda(z)
$$
for any $f_\lambda\in\CC(Z_0)^{(B_0)}_{\lambda}$. We can extend $\chi$ to a character $\lambda'\in M(T')$ of the torus $T'$. Put $F(t, z)=\lambda'(t)f_\lambda(z)$, $F\in\CC(T'\times Z_0)^{(T'\times B_0)}$. We have
$$
F(t k^{-1}, k z)=\lambda'(t)\lambda'(k)^{-1}\chi(k)f_\lambda(z)=F(t, z),
$$
so $F$ descends to a rational function on $Z$. And because $B\cap L = T' B_0$, the function $F$ is $(B\cap L)$-semi-invariant. We also have $F\vert_{Z_0}=f_\lambda$, hence the surjectivity of the map $\mathcal{M}_Z\to\mathcal{M}_{Z_0}$.

Now, if the character $\chi\in\mathcal{M}_Z$ of a $(B\cap L)$-eigenfunction maps to zero in $\mathcal{M}_{Z_0}$ then, as $B_0$ has an open orbit in $Z_0$, this function is constant on the fibers of $Z\to L/L_0$, and thus descends to an $L/L_0$-eigenfunction. The exactness in the middle term follows.

Thus we have the dual exact sequence
\begin{equation}
\label{eq: ses}
0\to\mathcal{N}_{Z_0}\to\mathcal{N}_{Z}\to N(L/L_0)\to 0.
\end{equation}
The valuation cone of $T'\times Z_0$ is $N(T')_\QQ\times\mathcal{V}_{Z_0}$, and, since $N(T')_\QQ\cong N(L/L_0)_\QQ$, the finite covering $T'\times Z_0\to Z$ gives $\mathcal{V}_Z=N(L/L_0)_\QQ \times \mathcal{V}_{Z_0}$ (under some splitting of~(\ref{eq: ses}) over $\QQ$). Everything else is straightforward.

\end{proof}

Finally, every $(B\cap L_0^\circ)$-stable prime divisor on $Z_0$ containing the closed orbit occurs as an irreducible component of $D\cap Z_0$ for a $(B\cap L)$-stable prime divisor $D$ on $Z$ containing the closed orbit. Each such component has the same valuation vector as $D$ under the embedding $\mathcal{N}_{Z_0}\hookrightarrow\mathcal{N}_Z$. Hence the cone $\mathcal{C}_{Z_0}$ spans $(\mathcal{N}_{Z_0})_\QQ\subset(\mathcal{N}_Z)_\QQ$ and coincides with $\mathcal{C}_Z=\mathcal{C}$.

\subsection{Star subdivision}
\label{s: star}

We describe a construction that generalizes the star subdivision of toric varieties (see \cite[\S\,3.3, ``Refinements of Fans and Blowups'']{CLS}).

Let $X$ be a simple spherical variety whose closed orbit $Y$ corresponds to a colored cone $(\mathcal{C},\mathcal{F})$. Suppose that $X$ is smooth and hence, in particular, locally factorial. There is a simple criterion for local factoriality of spherical varieties (see, e.g.,~\cite{Pas}).

\begin{proposition}
\label{prop: loc_factorial}
A simple spherical variety $X$ is locally factorial if and only if the colors in $\mathcal{F}$ have distinct images in $\mathcal{N}_X$, and there exists a basis $\{v_1,\dots,v_k\}\cup\rho(\mathcal{F})$ of the lattice $\mathcal{N}_X$ such that the cone $\mathcal{C}$ is generated by the family $\{v_1,\dots,v_l\}\cup\rho(\mathcal{F})$, where $l\leq k$.
\end{proposition}

Put $\mathcal{C}=\cone(u_1,\dots,u_n)$, where the $u_i\in\mathcal{N}$ are the primitive generators of the rays of $\mathcal{C}$. By Proposition~\ref{prop: loc_factorial}, these elements form a part of a basis of $\mathcal{N}$. Take a primitive vector $v\in\mathcal{C}\cap\mathcal{V}\cap\mathcal{N}$, and let $(\mathcal{C}_0,\mathcal{F}_0)$ be the colored face of $(\mathcal{C},\mathcal{F})$ such that $v\in\mathcal{C}_0^\circ$. We may assume that $\mathcal{C}_0=\cone(u_1,\dots,u_r)$. It may happen that $r=1$, $v=u_1$ and $\mathcal{C}_0=\cone(v)$ (see Example~\ref{exmpl: A2}). Then either $v=\rho(D)$ for some color $D\in\mathcal{F}$ and $\mathcal{F}_0=\{D\}$, or $v$ corresponds to a $G$-stable prime divisor of $X$ and $\mathcal{F}_0=\varnothing$. If $r>1$, we introduce the set of polyhedral cones
\begin{equation*}
F=\{\cone(A)\mid A\subset\{v,u_1,\dots,u_n\},\ \{u_1,\dots,u_r\}\not\subset A\}.
\end{equation*}

\begin{definition}
\label{def: star_cone}
The \textit{star subdivision} of a colored cone $(\mathcal{C},\mathcal{F})$ at a primitive vector $v\in\mathcal{C}\cap\mathcal{V}\cap\mathcal{N}$ is
\textup{(a)} the colored fan
\begin{equation*}
\Star((\mathcal{C},\mathcal{F}),v)=\{(\mathcal{C}',\rho^{-1}(\mathcal{C}')\cap\mathcal{F})\mid\mathcal{C}'\in F,\ (\mathcal{C}')^\circ\cap\mathcal{V}\neq\varnothing\},
\end{equation*}
if $v$ lies on no ray of $\mathcal{C}$, or \textup{(b)} the colored fan
\begin{equation*}
\Star((\mathcal{C},\mathcal{F}),v)=\{(\mathcal{C}',\mathcal{F}'\backslash\mathcal{F}_0)\mid(\mathcal{C}',\mathcal{F}')\in\FF_\mathcal{C}\},
\end{equation*}
where $\FF_\mathcal{C}$ is the fan of colored faces of $(\mathcal{C},\mathcal{F})$ (including $(\mathcal{C},\mathcal{F})$ itself), if $v$ spans a ray of $\mathcal{C}$.
\end{definition}

Now let $X$ be a spherical variety with the colored fan $\FF$, and let $v\in\Supp\FF\cap\mathcal{N}$ be a primitive vector. Suppose that $X$ is smooth along the center $Y$ of the valuation $v$, that is, at every point of $Y$. Then for every colored cone $(\mathcal{C},\mathcal{F})\in\FF$ with $v\in\mathcal{C}$ the star subdivision in the sense of Definition~\ref{def: star_cone} is well defined. Indeed, the subvariety $Y_0$ corresponding to $(\mathcal{C},\mathcal{F})$ is contained in the center $Y$ of $v$ (Proposition~\ref{prop: centers} and Remark~\ref{rm: inclusion_reversing}), so the simple subvariety of $X$ corresponding to $(\mathcal{C},\mathcal{F})$ is smooth, and Proposition~\ref{prop: loc_factorial} applies. This allows us to give the following definition.

\begin{definition}
\label{def: star_fan}
The \textit{star subdivision} of a colored fan $\FF$ at a primitive vector $v\in\Supp\FF\cap\mathcal{N}$ is the colored fan
\begin{equation*}
\Star(\FF,v)=\Biggl(\bigcup_{(\mathcal{C},\mathcal{F})\in\FF_v}\Star((\mathcal{C},\mathcal{F}),v)\Biggr)\cup(\FF~\backslash~\FF_v),
\end{equation*}
where $\FF_v=\{(\mathcal{C},\mathcal{F})\in\FF\mid v\in\mathcal{C}\}$.
\end{definition}

In other words, to construct the star subdivision of a colored fan one has to perform the regular star subdivision, but keep only the cones whose relative interior meets the valuation cone, and keep all the old colors unless a color falls on the subdividing ray.

\section{Main result}
\label{s: main_result}

The proof of the main result relies on Luna's étale slice theorem~\cite{Lu}: a smooth affine variety with an action of a reductive group $G$ and a unique closed orbit is isomorphic to a vector space with a linear action of $G$, provided that this closed orbit is a fixed point.

Let $V$ be a spherical $G$-module. It is a simple $G$-variety, and its fan is completely determined by the maximal colored cone $(\mathcal{C},\mathcal{F})$ corresponding to the zero orbit. Since $V$ is spherical, its simple submodules are pairwise non-isomorphic; denote them by $V_i$, so that $V=V_1\oplus\dots\oplus V_k$. We have the decomposition
\begin{equation}
\label{eq: grading}
\CC[V]=\bigoplus_{l_i\geq0}\CC[V]_{l_1,\dots,l_k}
\end{equation}
of $\CC[V]$ into the direct sum of subspaces consisting of the polynomials that are homogeneous of degree $l_i$ in the variables corresponding to $V_i$. Denote by $v_1,\dots,v_k$ the $G$-invariant valuations corresponding to this grading, that is, $v_i(f)=l_i$ for every nonzero $f\in\CC[V]_{l_1,\dots,l_k}$. Then we have the following lemma.

\begin{lemma}
\label{lemma: generators}
The vectors $v_1,\dots,v_k\in\mathcal{N}$ form the set of all primitive generators of the cone $\sigma=\mathcal{C}\cap l\mathcal{V}$ and a basis of the lattice $l\mathcal{V}\cap\mathcal{N}$.
\end{lemma}

\begin{proof}
The valuations $v$ in the linear part $l\mathcal{V}$ are distinguished among all invariant valuations by the property that the filtration $\{f\mid v(f)\geq c\}$, $c\in\QQ$, is in fact a grading of $\CC[V]$ by the value group $\ZZ\cong\Ima v\subset\QQ$. Namely, as $v\in l\mathcal{V}=\mathcal{T}_V^\perp$ (see Proposition~\ref{prop: valuation_cone_affine}), $v$ vanishes on all tails, and the grading is defined by $\deg f=\langle v,\lambda\rangle$ for nonzero $f\in\CC[V]_{(\lambda)}$. Among these valuations, those in $\mathcal{C}\cap l\mathcal{V}$ are the ones that have a center on $V$ (and define a nontrivial nonnegative grading), see Proposition~\ref{prop: centers}. Therefore the $v_i$ are primitive (since $v_i(f_{\lambda_i})=1$ for a highest weight vector $f_{\lambda_i}\in V_i^*$) and lie in $\mathcal{C}\cap l\mathcal{V}$.

Conversely, every $G$-invariant valuation $v$ takes a constant value $c_i$ on $V_i^*\backslash0\subset\CC[V]$, since the module $V_i^*$ is irreducible. If, moreover, $v$ lies in the linear part of the valuation cone and therefore defines a grading of $\CC[V]$, then $v$ is constant on each summand of the decomposition~\eqref{eq: grading}, so it is determined by its values on linear functions. Hence $v=c_1v_1+\dots+c_kv_k$. If $v$ is integral, that is, $v\in l\mathcal{V}\cap\mathcal{N}$, then $c_i\in\ZZ$ (being the value of $v$ on $V_i^*\backslash0$), and thus $v_1,\dots,v_k$ generate $l\mathcal{V}\cap\mathcal{N}$. If $v$ (not necessarily integral) has a center on $V$, that is, $v\in\mathcal{C}$, then $c_i\geq0$. Therefore $\sigma$ is contained in the cone generated by the vectors $v_i\in\sigma$, and so coincides with it. Finally, the $v_i$ are linearly independent. Hence they correspond bijectively to the rays of $\sigma$ and form a basis of $l\mathcal{V}\cap\mathcal{N}$.
\end{proof}

\begin{remark}
The proof of Lemma~\ref{lemma: generators} shows that the centers of the invariant valuations $v\in\sigma$ are the linear subspaces $U_1\oplus\dots\oplus U_k$, where $U_i$ is $0$ or $V_i$ depending on whether the vector $v_i$ does or does not belong to the set of generators of the smallest face of $\sigma$ containing $v$. For example, the center of $v_i$ is $\bigoplus_{j\neq i}V_j$, and the valuations in the relative interior $\sigma^\circ$ have the origin as their center.
\end{remark}

When the origin is blown up in a spherical module $V$, the exceptional divisor corresponds to the ray of $\FF(\Bl_0V)$ spanned by the vector $v_0$ of the valuation that assigns to a polynomial $f$ the degree of its lowest homogeneous component with respect to the total degree. That is, $v_0(\CC[V]_{l_1,\dots,l_k}\backslash0)=l_1+\dots+l_k$. By Lemma~\ref{lemma: generators}, $v_0$ is the sum $v_1+\dots+v_k$ of the primitive vectors on the rays of $\mathcal{C}\cap l\mathcal{V}$.

We are now ready to state and prove the main result.

\begin{theorem}
\label{thm: main}
Let $Y$ be a smooth $G$-stable irreducible subvariety of a spherical $G$-variety $X$, and let $(\mathcal{C},\mathcal{F})$ be the corresponding colored cone of $\FF(X)$. Suppose that $X$ is smooth at the points of $Y$. Put
\begin{align*}
S_\mathcal{F}&=\{\alpha\in S\mid\mathcal{D}(\alpha)\subset\mathcal{F}\},\\
\Sigma_*&=\Sigma\cap\langle S_\mathcal{F}\rangle_\QQ,\\
\mathcal{V}_*&=-\Sigma_*^\vee.
\end{align*}
Let $v_0\in\mathcal{N}$ be the sum of the primitive generators of the rays of the cone $\mathcal{C}\cap l\mathcal{V}_*$. Then for the blowup $X'=\Bl_YX$ we have $\FF(X')=\Star(\FF(X),v_0)$.
\end{theorem}

\begin{proof}
Since blowing up is a local operation, we may replace $X$ by the open $G$-stable subset
\begin{equation*}
X_0=\bigcup_{Y'\in\mathcal{Y}}X_{Y'}
\end{equation*}
containing $Y$, where $\mathcal{Y}$ is the set of closed $G$-orbits in $Y$, and thus assume that $X$ is smooth. The fan $\FF(X)$ is then replaced by the smaller colored fan $\FF(X_0)$, whose maximal colored cones correspond to the orbits in $\mathcal{Y}$. So from now on $X$ is assumed to be smooth. Note that all colored cones of $\FF(X)$ corresponding to $G$-stable subvarieties of $X$ not contained in $Y$ remain in $\FF(X')$ unchanged.

As already noted, the exceptional divisor $E$ of the blowup of a smooth variety along a smooth center is prime. We show that its valuation vector $v_E$ equals $v_0$.

Apply the local structure theorem to $X\supset Y$ (see \S\,\ref{s: local_structure}). We obtain an $L$-spherical variety $Z$ and an $L_0$-stable variety $Z_0$ with a fixed point, spherical with respect to the identity component $L_0^\circ$ of $L_0$. Both are smooth, since $X$ is smooth. By Luna's étale slice theorem~\cite{Lu}, $Z_0$ is equivariantly isomorphic to a spherical $L_0^\circ$-module $V$. Hence, in order to compute $v_E$, the problem of blowing up $X$ along $Y$ reduces to blowing up $X^\circ$ along $Y^\circ$, then, since $Y^\circ=P_u\times(Y\cap Z)\subset X^\circ = P_u\times Z$, to blowing up $Z$ along $Y\cap Z$, and, finally, to blowing up the origin in the spherical module $V$. Namely, if $v$ is the $L_0^\circ$-invariant valuation corresponding to the exceptional divisor of the blowup of $0\in V$, then $v_E=v\in l\mathcal{V}_V=l\mathcal{V}_{Z_0}\subset l\mathcal{V}_Z\subset\mathcal{V}_Z$ (Lemma~\ref{lemma: localization}). The equality $v=v_0$ follows from Lemma~\ref{lemma: generators}, because, by \S\,\ref{s: local_structure}, the cone $\mathcal{C}_V=\mathcal{C}_{Z_0}$ of the zero orbit in $V$ is equal to the cone $\mathcal{C}$ of $Y$, we have the equality $\mathcal{V}_Z=\mathcal{V}_*$, and the embedding $\mathcal{N}_{Z_0}\subset\mathcal{N}_Z$ is saturated (by Lemma~\ref{lemma: localization}). Thus, $v_E=v_0$.

It remains to find out what happens to a colored cone $(\mathcal{C}(Y_0),\mathcal{F}(Y_0))\in\FF(X)$ with $Y_0\subset Y$. Recall that $\pi\colon X'\to X$ is the blowup morphism. Let $Y'_0$ be an arbitrary irreducible closed $G$-stable subvariety of $\pi^{-1}(Y_0)\subset E$ that maps onto $Y_0$. Then $\mathcal{C}(Y'_0)$ is generated by $v_E$ and some rays of~$\mathcal{C}(Y_0)$.

Suppose first that $v_E$ does not lie on any ray of $\mathcal{C}(Y_0)$. Then $\mathcal{F}(Y'_0)$ consists of those colors $D\in\mathcal{F}(Y_0)$ whose images lie in $\mathcal{C}(Y'_0)$. Indeed, suppose that there is a color $D\in\mathcal{F}(Y_0)$ such that $Y'_0$ is not contained in $\overline{D}\subset X'$ (the strict transform of the closure of $D$ in $X$). Since $X'$ is smooth, Proposition~\ref{prop: loc_factorial} implies that the corresponding ray of $\mathcal{C}(Y'_0)$ is spanned by the valuation $v$ of some $G$-stable prime divisor of $X'$. But the $G$-stable prime divisors of $X'$ are the strict transforms of the $G$-stable prime divisors of $X$, together with $E$. By assumption $v_E\neq v$, and the valuation of a $G$-stable prime divisor of $X$ cannot span a ray of $(\mathcal{C}(Y_0),\mathcal{F}(Y_0))$ that contains the image of a color from $\mathcal{F}(Y_0)$ (this follows, for instance, from the smoothness of $X$ and Proposition~\ref{prop: loc_factorial}, although it holds for any spherical variety $X$, smooth or not, and any irreducible $G$-stable subvariety $Y$). This is a contradiction.

Suppose now that $v_E$ spans a ray of $\mathcal{C}(Y_0)$. This ray cannot be spanned by the valuation of a $G$-stable prime divisor on $X$, because distinct $G$-invariant valuations have distinct images in $\mathcal{N}_\QQ$. Therefore the image of some color $D\in\mathcal{F}(Y_0)$ must lie on it. By Proposition~\ref{prop: loc_factorial}, this color does not belong to $\mathcal{F}(Y'_0)$, and applying the previous argument to the colors on the remaining rays, we obtain
\begin{equation*}
\mathcal{F}(Y'_0)=\{D'\in\mathcal{F}(Y_0)\mid\rho(D')\in\mathcal{C}(Y'_0)\}\backslash\{D\}.
\end{equation*}

Thus, in either case, the colored cone of $Y'_0$ lies in $\Star(\FF(X),v_0)$. Conversely, since the morphism $\pi\colon X'\to X$ is proper, the supports of $\FF(X)$ and $\FF(X')$ coincide (Corollary~\ref{cor: proper}), and hence every colored cone of $\Star(\FF(X),v_0)$ must belong to $\FF(X')$.
\end{proof}

\section{Examples}
\label{s: examples}

\begin{figure}[ht]
\centering

\pgfmathsetmacro{\a}{1.7}
\pgfmathsetmacro{\ang}{60}
\pgfmathsetmacro{\h}{\a*sin(\ang)}

\begin{tikzpicture}
 
\begin{scope}
    \path[clip] (-2.7*\a,-1.8*\h) rectangle (2.1*\a,2.5*\h);
    
    \foreach \i in {-3,...,3} {
        \draw[style=help lines,dashed] (\i*\a,0) -- ++(\ang:3*\a);
        \draw[style=help lines,dashed] (\i*\a,0) -- ++(4*\ang:3*\a);
        \draw[style=help lines,dashed] (\i*\a,0) -- ++(-\ang:3*\a);
        \draw[style=help lines,dashed] (\i*\a,0) -- ++(2*\ang:3*\a);
    }
    \foreach \j in {-2,...,2} {
        \draw[style=help lines,dashed] (0,\j*\h) -- ++(0:3*\a);
        \draw[style=help lines,dashed] (0,\j*\h) -- ++(3*\ang:3*\a);
    }

    \foreach \i in {-3,...,3} {
        \foreach \j in {-3,...,3} {
            \node[draw,circle,inner sep=1pt,fill=gray,gray] at (\i*\a+\j*\a/2,\j*\h) {};
        }
    }

    \node[above right,gray] at (1.1*\a,0) {$2\varepsilon_1$};
    \node[above left,gray] at (-0.5*\a,\h) {$2\varepsilon_2$};
    \node[above left,gray] at (-0.5*\a,-\h) {$2\varepsilon_3$};

    \draw[thin,-latex,blue, fill=blue, fill opacity=0.1] (0,0) -- (0,3*\h) -- (-3*\a,3*\h) -- (-3*\a,-3*\a) -- (-1.5*\a,-3*\h) -- cycle;

    \draw[blue,ultra thick,opacity=0.5] (0,0.5*\a) arc (1.5*\ang:4*\ang:0.5*\a);

    \node[blue] at (-2*\a,1.5*\h) {$\mathcal{C}$};
    
    \draw[thin,-latex,gray, fill=gray, fill opacity=0.3] (0,0) -- (-3*\a,0) -- (-3*\a,-3*\a) -- (-1.5*\a,-3*\h) -- cycle;

    \draw[gray,ultra thick] (-0.75*\a,0) arc (3*\ang:4*\ang:0.75*\a);

    \node[above,gray] at (-2*\a,-1.5*\h) {$\mathcal{V}$};

    \draw[ultra thick,-latex,black] (0,0) -- (1.5*\a,-\h) node[pos=0.5, below] {$2\alpha_1$};

    \node[draw,circle,inner sep=1.5pt,fill=orange,orange] at (1.5*\a,-\h) {};
    \node[right] at (1.5*\a,-\h) {$D_1$};

    \draw[ultra thick,-latex] (0,0) -- (0,2*\h) node[pos=0.5, below right] {$2\alpha_2$};

    \node[draw,circle,inner sep=1.5pt,fill=orange,orange] at (0,2*\h) {};
    \node[right] at (0,2*\h) {$D_2$};


    \node[draw,circle,inner sep=1.5pt,fill] at (-0.5*\a,-\h) {};
    \node[right] at (-0.5*\a,-\h) {$D_3$};

    \draw[ultra thick,-latex] (0,0) -- (-\a,0) node[pos=0.7,above] {$v$};
    \node[draw,circle,inner sep=1.5pt,fill] at (-\a,0) {};
    \node[above left] at (-\a,0) {$E$};
    
\end{scope}
\end{tikzpicture}
\caption{The variety of conics}
\label{fig: conics}
\end{figure}
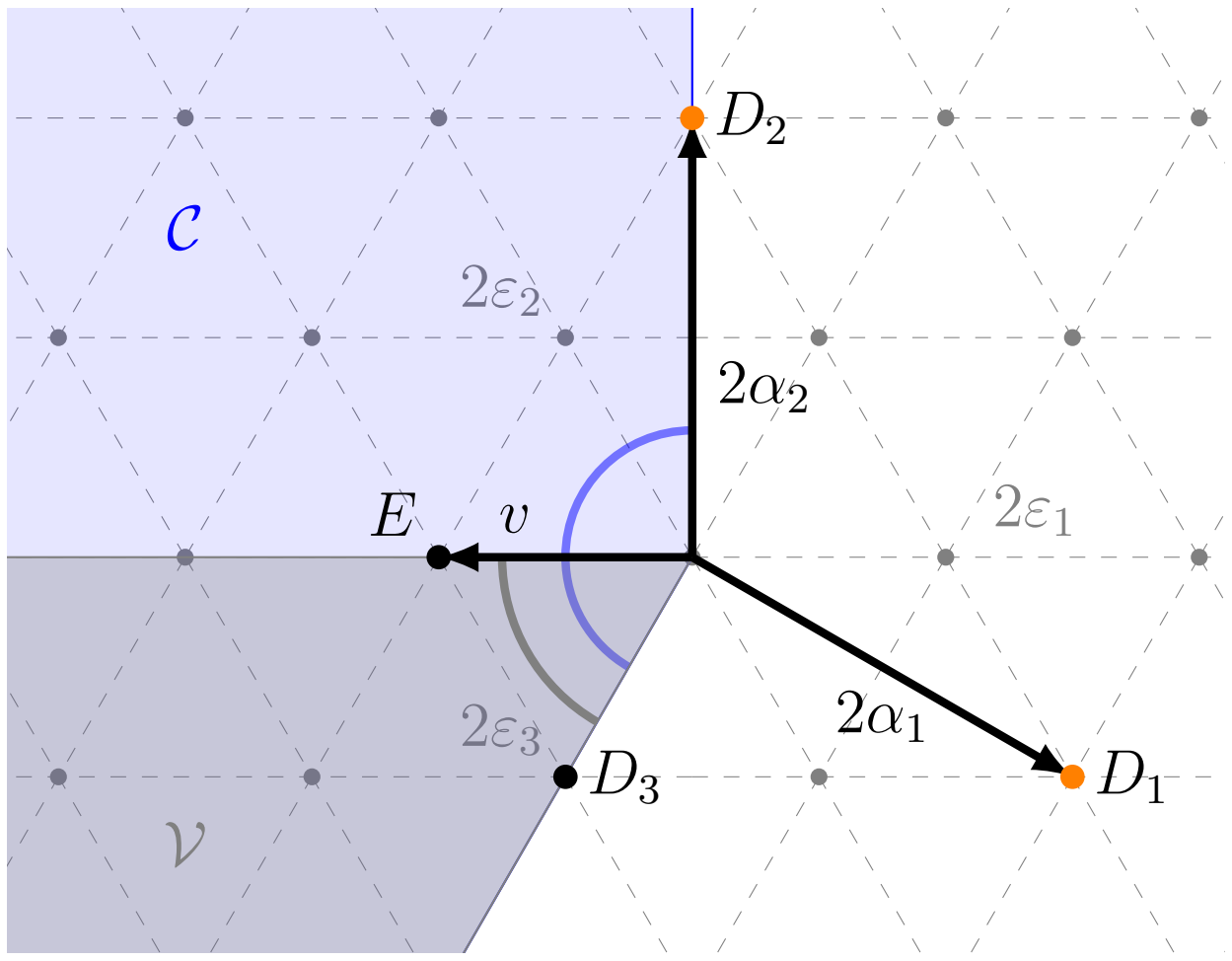

\begin{example}
\label{exmpl: conics}
Consider the space of smooth conics in $\PP^2$. It is a spherical homogeneous space $G/H$ for the group $G=\GL_3$. Fix the Borel subgroup $B$ of upper triangular matrices and the maximal torus $T$ of diagonal matrices. We may regard $G/H$ as an open subset of $X=\PP^5$, the projectivization of the space of symmetric $3\times3$ matrices. Let $\Delta_1=x_{11}$, $\Delta_2=x_{11}x_{22}-x_{12}^2$ and $\Delta_3=\det x$ be the leading principal minors of a symmetric matrix $x=(x_{ij})$. The first two are $B$-semi-invariant of weights $2\eps_1$ and $2(\eps_1+\eps_2)$, and $\Delta_3$ is $G$-semi-invariant of weight $2(\eps_1+\eps_2+\eps_3)$. Here we identify the character lattices of $B$ and $T$ and put $\eps_i(t)=t_i$ for a diagonal matrix $t=\diag(t_1,t_2,t_3)$. The equations $\Delta_i=0$ define divisors $D_1$, $D_2$ and $D_3$ in $\PP^5$; the divisors $D_1$ and $D_2$ correspond to colors of $G/H$, while $D_3$ is the only divisor in the complement of the open $G$-orbit.

The lattice $\mathcal{M}\cong\ZZ^2$ consists of the weights $\lambda=2(l_1\eps_1+l_2\eps_2+l_3\eps_3)$ of the rational functions $f_\lambda=\Delta_1^{l_1-l_2}\Delta_2^{l_2-l_3}\Delta_3^{l_3}\in\CC(\PP^5)$ with $l_1+l_2+l_3=0$. One can show that there are two spherical roots, $2\alpha_1$ and $2\alpha_2$, where $\alpha_i=\eps_i-\eps_{i+1}$ are the simple roots of $\GL_3$ (see, e.g., \cite[Example 17.12]{Ti}). A direct computation shows that $v_{D_i}=(\eps_i^*-\eps_{i+1}^*)/2$ for $i=1,2$, that $v_{D_3}=-(v_{D_1}+2v_{D_2})/3$, and that $\mathcal{D}(\alpha_1)=\{D_1\}$ and $\mathcal{D}(\alpha_2)=\{D_2\}$. Here $\eps_i^*$ denotes the restriction to the subspace $\mathcal{M}_\QQ$ of the linear function dual to $\eps_i$ in the basis of the character lattice of $T$. The valuation cone $\mathcal{V}$ is generated by two rays: one is spanned by $v_{D_3}$, and the other by the vector $v:=-(2v_{D_1}+v_{D_2})/3$.

The variety $X$ has a unique closed $G$-orbit $Y$, namely, the image in $\PP^5$ of the set of matrices of rank $1$; it lies in $D_2$ and $D_3$. As a simple spherical variety, $X$ is determined by its maximal colored cone $(\mathcal{C},\mathcal{F})$, where $\mathcal{C}=\cone(v_{D_2},v_{D_3})$ and $\mathcal{F}=\{D_2\}$. All these objects are shown in Figure~\ref{fig: conics} (the lattices $\mathcal{M}$ and $\mathcal{N}$ are superimposed in one picture in a suitable way).

The subvariety $Y$ is smooth, and in the notation of Theorem~\ref{thm: main} we have $S_\mathcal{F}=\{\alpha_2\}$ and $\Sigma_*=\{2\alpha_2\}$, so $\mathcal{C}\cap l\mathcal{V}_*=\cone(v)$. Hence the exceptional divisor $E$ of the blowup $X'$ of $X$ along $Y$ corresponds to the valuation $v_0=v$, and the colored fan $\FF(X')$ equals $\Star((\mathcal{C},\mathcal{F}),v)=\{(\mathcal{V},\varnothing),(\cone(v_{D_3}),\varnothing),(\cone(v),\varnothing),(0,\varnothing)\}$.
\end{example}

\begin{example}
\label{exmpl: A2}
It may happen that the image $\rho(D)$ of a color $D\in\mathcal{F}$ of a colored cone $(\mathcal{C},\mathcal{F})$ of $\FF(X)$ lies on a new ray of the fan $\FF(X')$ of the blowup. In this case the strict transform of $D$ in $X'$ contains no $G$-orbit, and $D$ does not occur among the colors of the colored cones of the new fan.

Consider the affine plane $X=\AAA^2$ with the action of $G=\SL_2$. It has a unique $B$-stable divisor $D$, given by the polynomial $f(x,y)=y$ of weight $\eps_1$ with respect to the Borel subgroup of upper triangular matrices; here $\eps_1(t)=t_1$ for $t=\diag(t_1,t_1^{-1})\in\SL_2$. Then $\mathcal{M}=\ZZ\eps_1$ and $\mathcal{D}=\{D\}$, and one can show that $\mathcal{V}=\mathcal{N}_\QQ$ (for example, note that the open $G$-orbit $\AAA^2\backslash0$ is of the form $G/U$, where $U\subset B$ is a maximal unipotent subgroup, so $X$ is a so-called \textit{horospherical} variety \cite[\S\,2.7]{Ti}, and the set of spherical roots of such a variety is empty \cite[Theorem 21.10]{Ti}). The variety $X$ has a unique closed orbit $Y=\{0\}$, and $\FF(X)=\{(\cone(v_D),\{D\}),(0,\varnothing)\}$, where $\pairing{v_D}{\eps_1}=1$.

\begin{figure}[ht]
\centering
\pgfmathsetmacro{\a}{1.4}
\begin{subfigure}{0.5\textwidth}
\centering
\begin{tikzpicture}
\begin{scope}

    \draw[style=help lines,dashed] (-2.2*\a,0) -- (2.2*\a,0);
    
    \draw[] (0,-0.07*\a) -- (0,0.07*\a);

    \draw[] (0,0) -- (2.2*\a,0);
    
    \node[draw,circle,inner sep=1.5pt,fill=orange,orange] at (\a,0) {};

    \node[draw,circle,inner sep=1pt,fill=gray,gray] at (-\a,0) {};

    \node[draw,circle,inner sep=1pt,fill=gray,gray] at (-2*\a,0) {};

    \node[draw,circle,inner sep=1pt,fill=gray,gray] at (2*\a,0) {};

    \node[above right] at (\a,0) {$D$};
    \node[below,gray] at (\a,0) {$\varepsilon_1$};
    \node[below,gray] at (0,0) {$0$};
    
\end{scope}
\end{tikzpicture}
\subcaption{$\mathfrak{F}(\AAA^2)$}
\end{subfigure}%
\begin{subfigure}{0.5\textwidth}
\centering
\begin{tikzpicture}

\begin{scope}

    \draw[style=help lines,dashed] (-2.2*\a,0) -- (2.2*\a,0);
    
    \draw[] (0,-0.07*\a) -- (0,0.07*\a);

    \draw[] (0,0) -- (2.2*\a,0);
    
    \node[draw,circle,inner sep=1.5pt,fill] at (\a,0) {};

    \node[draw,circle,inner sep=1pt,fill=gray,gray] at (-\a,0) {};

    \node[draw,circle,inner sep=1pt,fill=gray,gray] at (-2*\a,0) {};

    \node[draw,circle,inner sep=1pt,fill=gray,gray] at (2*\a,0) {};

    \node[above right] at (\a,0) {$E$};
    \node[below,gray] at (\a,0) {$\varepsilon_1$};
    \node[below,gray] at (0,0) {$0$};

    \draw[thick,-latex,black] (0,0) -- (\a,0) node[pos=0.5, above] {$v$};
    
\end{scope}
\end{tikzpicture}
\subcaption{$\mathfrak{F}(\Bl_0 \AAA^2)$}
\end{subfigure}
\caption{Blowup of an $\SL_2$-spherical affine plane}
\label{fig: affine_plane}
\end{figure}
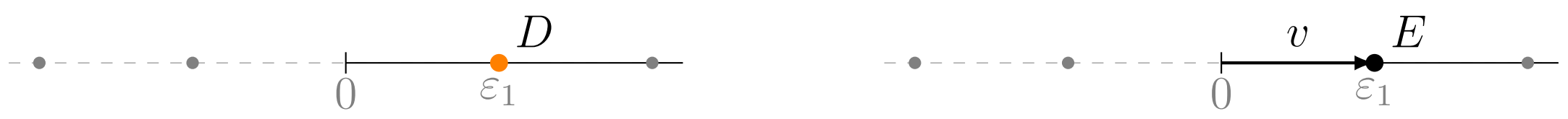

Since the set $\Sigma$ of spherical roots is empty, for the blowup $X'=\Bl_0\AAA^2$ of the closed orbit corresponding to the colored cone $(\cone(v_D),\{D\})$ we have $l\mathcal{V}_*=\mathcal{V}_*=\mathcal{V}$, and the exceptional divisor $E$ corresponds to the valuation $v$ whose image in $\mathcal{N}$ coincides with $\rho(D)$. Moreover, $\FF(X')=\Star(\FF(X),v)=\{(\cone(v),\varnothing),(0,\varnothing)\}$. See Figure~\ref{fig: affine_plane}.
\end{example}

\begin{example}
\label{exmpl: wedge}
Consider the $G$-module $V=\bigwedge^2(\CC^4)^*\oplus(\CC^4)^*$, where $G=\GL_4\times\CC^\times$; the group $\GL_4$ acts on both summands, and $\CC^\times$ acts on $\CC^4$ by dilations. In $\GL_4$ we fix the Borel subgroup of upper triangular matrices and the maximal torus of diagonal matrices. Their products with $\CC^\times$ are a Borel subgroup and a maximal torus of $G$, respectively. The module $V$ is spherical, and the lattice $\mathcal{M}_V$ of its $B$-eigenfunctions is generated by the weights $\omega_1+\eps$, $\omega_2$, $\omega_3+\eps$ and $\omega_4$ of regular semi-invariant polynomials of degrees $1,1,2$ and $2$, respectively (see the list in~\cite{Kn_MFS}). Here $\eps$ is the basic character of $\CC^\times$, and $\omega_i=\eps_1+\dots+\eps_i$ is the highest weight of the $\GL_4$-module $\bigwedge^i\CC^4$, $i=1,\dots,4$; we identify the characters of $B$ and $T$ and put $\eps_i(t)=t_i$ for $t=\diag(t_1,\dots,t_4)$. There are two spherical roots, $\beta_1=\eps_2-\eps_4=\alpha_2+\alpha_3$ and $\beta_2=\eps_1-\eps_3=\alpha_1+\alpha_2$, where $\alpha_i=\eps_i-\eps_{i+1}$ are the simple roots of $\GL_4$ (see~\cite{Kn_MFS}).

\begin{figure}[ht]
\centering
\begin{subfigure}{0.5\textwidth}
\centering
\begin{tikzpicture}[scale=0.3]
    \path[clip] (-10cm,-10cm) rectangle (10cm,10cm);
    \coordinate (O) at (0,0);
    \coordinate (e1) at (-3,-1.5);
    \coordinate (e2) at (2,-2);
    \coordinate (e3) at (0,2);
    
    \coordinate (v1) at ($-1*(e1)$);
    \coordinate (v2) at ($(e1)+(e2)-2*(e3)$);
    \coordinate (v3) at ($(e1)-(e2)+(e3)$);
    
    \coordinate (C1) at ($2.5*(v1)$);
    \coordinate (C2) at ($1.1*(v2)$);
    \coordinate (C3) at ($2*(v3)$);

    \coordinate (D12) at ($0.9*(C1)+0.9*(C2)$);
    \coordinate (D23) at ($0.9*(C2)+0.7*(C3)$);
    \coordinate (D31) at ($0.8*(C3)+0.9*(C1)$);

    \coordinate (int13) at ($-1*(e2)+1*(e3)$);
    \coordinate (I1) at ($2*(int13)$);

    \coordinate (G1) at (intersection of  C3--D31 and O--I1);

    \coordinate (E0) at ($2.5*(e3)$);
    \coordinate (H1) at ($-2*(e2)$);
    \coordinate (E1) at ($(E0)+(H1)$);

    \coordinate (J1) at (intersection of E0--E1 and C3--D31);
    \coordinate (K1) at (intersection of E1--H1 and C3--D31);

    \coordinate (H2) at ($-3*(e1)$);

    \coordinate (E2) at ($(H2)+(E0)$);

    \coordinate (J2) at (intersection of E0--E2 and C1--D31);

    \coordinate (F0) at ($-1*(E0)$);
    \coordinate (F1) at ($(F0)+(H1)$);
    \coordinate (F2) at ($(H2)+(F0)$);

    \coordinate (G2) at (intersection of C1--D12 and F0--F2);
    
    \begin{scope}
    \foreach \i in {-3,...,3} {
        \foreach \j in {0,...,3} {
            \foreach \k in {0} {
                \node[draw,circle,inner sep=0.1pt,fill] at ($\i*(e1)+\j*(e2)+\k*(e3)$) {};
            }
        }
    }

    \foreach \i in {0,...,3} {
        \foreach \j in {-3,...,0} {
            \foreach \k in {0} {
                \node[draw,circle,inner sep=0.1pt,fill] at ($\i*(e1)+\j*(e2)+\k*(e3)$) {};
            }
        }
    }

    \foreach \i in {-3,...,-1} {
        \foreach \j in {-3,...,-1} {
            \foreach \k in {0} {
                \node[circle,inner sep=0.1pt,fill=gray,opacity=0.5] at ($\i*(e1)+\j*(e2)+\k*(e3)$) {};
            }
        }
    }

    \foreach \i in {-3,...,-2} {
        \foreach \j in {-3,...,-2} {
            \foreach \k in {0} {
                \node[draw,circle,inner sep=0.1pt,fill] at ($\i*(e1)+\j*(e2)+\k*(e3)$) {};
            }
        }
    }

    \draw[style=help lines] ($3*(e1)$) -- ($-3*(e1)$);
    \draw[style=help lines] ($3*(e2)$) -- ($-3*(e2)$);
    \draw[style=help lines] ($3*(e3)$) -- ($-3*(e3)$);

    \fill[blue!50,fill opacity=0.15] (O) -- (C1) -- (D12) -- (C2) -- cycle;
    \draw[thick,blue,opacity=0.5] (O) -- (C2);
    \fill[blue,fill opacity=0.15] (O) -- (C2) -- (D23) -- (C3) -- cycle;
    \draw[thick,blue,opacity=0.5] (O) -- (C3);
    
    \draw[thick,blue,opacity=0.5] (O) -- (C1);

    \draw[thick,blue,opacity=0.2] (O) -- (G1);
    \fill[blue!70,fill opacity=0.15] (O) -- (C3) -- (G1) -- cycle;
    \fill[blue!70,fill opacity=0.15] (E0) -- (J1) -- (D31) -- (J2) -- cycle;

    \fill[gray!150,fill opacity=0.2] (O) -- (E0) -- (E1) -- (K1) -- (G1) -- cycle;
    \fill[gray!150,fill opacity=0.05] (O) -- (G1) -- (K1) -- (H1) -- cycle;
    \fill[gray,fill opacity=0.2] (O) -- (E0) -- (E2) -- (H2) -- cycle;
    \draw[thick,gray,opacity=0.5] (O) -- (E0);

    \fill[blue!70,fill opacity=0.05] (O) -- (G1) -- (J1) -- (E0) -- (J2) -- (C1) -- cycle;

    \fill[gray!150,fill opacity=0.05] (O) -- (F0) -- (F1) -- (H1) -- cycle;
    \fill[gray,fill opacity=0.05] (O) -- (F0) -- (G2) -- (C1) -- cycle;
    \fill[gray,fill opacity=0.2] (C1) -- (G2) -- (F2) -- (H2) -- cycle;
    \draw[thick,gray,opacity=0.5] (O) -- (F0);

    \draw[thick,-latex,black] (O) -- (e1) node[above] {$\beta_1$};
    \draw[thick,-latex,black] (O) -- (e2) node[pos=0.8, right] {$\beta_2$};
    \draw[thick,-latex,black] (O) -- (e3) node[right] {$\gamma$};

    \node[draw,circle,inner sep=1.5pt,fill=orange,orange] at (v2) {};
    \node[right] at (v2) {$D_2$};

    \node[draw,circle,inner sep=1.5pt,fill=orange,orange] at (v3) {};
    \node[above] at (v3) {$D_3$};

    \node[draw,circle,inner sep=1.5pt,fill=black,black] at (v1) {};
    \node[above] at (v1) {$D_4$};

    \node[right,gray] at ($-2*(e1)+(e3)$) {$\mathcal{V}$};
    \node[below,blue] at ($2*(e1)$) {$\mathcal{C}_1$};

    \end{scope}
\end{tikzpicture}
\subcaption{$\mathfrak{F}(X_{Y_1})$}
\end{subfigure}%
\begin{subfigure}{0.5\textwidth}
\centering
\begin{tikzpicture}[scale=0.3]
    \path[clip] (-10cm,-10cm) rectangle (10cm,10cm);
    \coordinate (O) at (0,0);
    \coordinate (e1) at (-3,-1.5);
    \coordinate (e2) at (2,-2);
    \coordinate (e3) at (0,2);
    
    \coordinate (v1) at ($-1*(e1)$);
    \coordinate (v2) at ($(e1)+(e2)-2*(e3)$);
    \coordinate (v3) at ($(e1)-(e2)+(e3)$);
    
    \coordinate (C1) at ($2.5*(v1)$);
    \coordinate (C2) at ($1.1*(v2)$);
    \coordinate (C3) at ($2*(v3)$);

    \coordinate (D12) at ($0.9*(C1)+0.9*(C2)$);
    \coordinate (D23) at ($0.9*(C2)+0.7*(C3)$);
    \coordinate (D31) at ($0.8*(C3)+0.9*(C1)$);

    \coordinate (int13) at ($-1*(e2)+1*(e3)$);
    \coordinate (I1) at ($2*(int13)$);

    \coordinate (G1) at (intersection of  C3--D31 and O--I1);

    \coordinate (E0) at ($2.5*(e3)$);
    \coordinate (H1) at ($-2*(e2)$);
    \coordinate (E1) at ($(E0)+(H1)$);

    \coordinate (J1) at (intersection of E0--E1 and C3--D31);
    \coordinate (K1) at (intersection of E1--H1 and C3--D31);

    \coordinate (H2) at ($-3*(e1)$);

    \coordinate (E2) at ($(H2)+(E0)$);

    \coordinate (J2) at (intersection of E0--E2 and C1--D31);

    \coordinate (F0) at ($-1*(E0)$);
    \coordinate (F1) at ($(F0)+(H1)$);
    \coordinate (F2) at ($(H2)+(F0)$);

    \coordinate (G2) at (intersection of C1--D12 and F0--F2);

    \coordinate (H3) at ($-1*(H1)$);
    \coordinate (E3) at ($(E0)+(H3)$);
    \coordinate (F3) at ($(H3)-(E0)$);

    \coordinate (int12) at ($1*(e2)-2*(e3)$);
    \coordinate (I3) at ($2*(int12)$);
    \coordinate (G3) at (intersection of C2--D12 and O--I3);
    \coordinate (K3) at (intersection of F3--F0 and C2--D12);
    \coordinate (S3) at (intersection of E3--F3 and C2--D12);
    \coordinate (T3) at (intersection of O--C1 and E0--E3);
    
    \begin{scope}
    \foreach \i in {-3,...,3} {
        \foreach \j in {0,...,3} {
            \foreach \k in {0} {
                \node[draw,circle,inner sep=0.1pt,fill] at ($\i*(e1)+\j*(e2)+\k*(e3)$) {};
            }
        }
    }

    \foreach \i in {0,...,3} {
        \foreach \j in {-3,...,0} {
            \foreach \k in {0} {
                \node[draw,circle,inner sep=0.1pt,fill] at ($\i*(e1)+\j*(e2)+\k*(e3)$) {};
            }
        }
    }

    \foreach \i in {-1} {
        \foreach \j in {-3,...,-1} {
            \foreach \k in {0} {
                \node[circle,inner sep=0.1pt,fill=gray,opacity=0.5] at ($\i*(e1)+\j*(e2)+\k*(e3)$) {};
            }
        }
    }

    \foreach \i in {-3,...,-2} {
        \foreach \j in {-3,...,-1} {
            \foreach \k in {0} {
                \node[draw,circle,inner sep=0.1pt,fill] at ($\i*(e1)+\j*(e2)+\k*(e3)$) {};
            }
        }
    }

    \draw[style=help lines] ($3*(e1)$) -- ($-3*(e1)$);
    \draw[style=help lines] ($3*(e2)$) -- ($-3*(e2)$);
    \draw[style=help lines] ($3*(e3)$) -- ($-3*(e3)$);

    \fill[blue!50,fill opacity=0.15] (O) -- (C2) -- (G3) -- cycle;
    \fill[blue!50,fill opacity=0.15] (T3) -- (C1) -- (D12) -- (S3) -- (E3) -- cycle;
    \draw[thick,blue,opacity=0.5] (O) -- (C2);
    \fill[blue,fill opacity=0.15] (O) -- (C2) -- (D23) -- (C3) -- cycle;
    \draw[thick,blue,opacity=0.5] (O) -- (C3);
    
    \draw[thick,blue,opacity=0.5] (O) -- (C1);

    \draw[thick,blue,opacity=0.2] (O) -- (G1);
    \fill[blue!70,fill opacity=0.15] (O) -- (C3) -- (G1) -- cycle;
    \fill[blue!70,fill opacity=0.15] (T3) -- (J1) -- (D31) -- (C1) -- cycle;
    \fill[blue!70,fill opacity=0.05] (O) -- (G1) -- (J1) -- (T3) -- cycle;

    \fill[gray!150,fill opacity=0.2] (O) -- (E0) -- (E1) -- (K1) -- (G1) -- cycle;
    \fill[gray!150,fill opacity=0.05] (O) -- (G1) -- (K1) -- (F1) -- (K3) -- (G3) -- cycle;

    \fill[gray!150,fill opacity=0.2] (O) -- (G3) -- (K3) -- (F3) -- (E3) -- (E0) -- cycle;

    \fill[blue!50,fill opacity=0.05] (O) -- (G3) -- (S3) -- (E3) -- (T3) -- cycle;

    \draw[thick,-latex,black] (O) -- (e1) node[above] {$\beta_1$};
    \draw[thick,-latex,gray] (O) -- ($(e2)-2*(e3)$) node[pos=0.8, right] {$v_1$};
    \draw[thick,-latex,gray] (O) -- ($(e3)-(e2)$) node[pos=1, right] {$v_2$};

    \draw[thick,-latex,black] (O) -- ($-1*(e3)$) {};

    \node[] at ($-1.2*(e3)+0.2*(e2)$) {$v_0$};



    \draw[thick,blue,opacity=0.2] (O) -- (G3);

    \node[above right,gray] at (e2) {$\mathcal{V}_*$};
    \node[below,blue] at ($2*(e1)$) {$\mathcal{C}_1$};

    \end{scope}
\end{tikzpicture}
\subcaption{Localization of $\mathfrak{F}(X_{Y_1})$ at the cone $(\mathcal{C}_1,\mathcal{F}_1)$}
\end{subfigure}
\caption{Blowup of the closed orbit $Y_1$ in $\PP(\bigwedge^2(\CC^4)^*\oplus(\CC^4)^*)$}
\label{fig: interesting_example}
\end{figure}
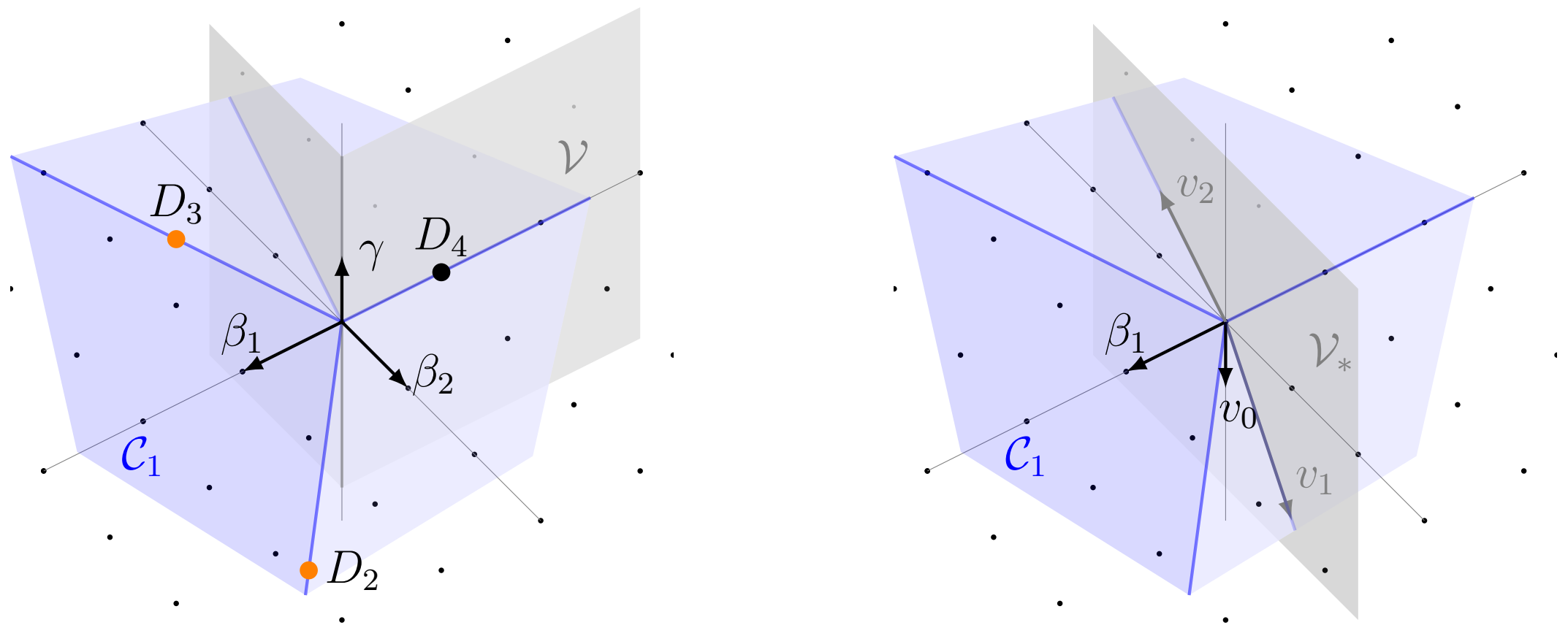

We now consider the projective spherical $G$-variety $X=\PP(V)$. Its weight lattice $\mathcal{M}=\mathcal{M}_X$ is the kernel of the map $\mathcal{M}_V\to\ZZ$ that sends the weight $\lambda$ of a regular function $f_\lambda$ to the degree of $f_\lambda$ (extended by linearity). One checks that $\mathcal{M}$ is generated by the spherical roots $\beta_1=\omega_2-(\omega_1+\eps)-\omega_4+(\omega_3+\eps)$ and $\beta_2=(\omega_1+\eps)+\omega_2-(\omega_3+\eps)$ and the vector $\gamma=(\omega_3+\eps)-2\omega_2$ (all of them lie in $\mathcal{M}$, and together with $\omega_2$ they form a basis of $\mathcal{M}_V$). Let $e_i$ be the standard basis vectors of $\CC^4$, and $e_i^*$ be the dual basis. Set $x=\sum_{i<j} x_{ij} e_i^*\wedge e_j^*$ and $y=\sum_i y_i e_i^*$. The equations with the left-hand sides
\begin{gather*}
f_{\omega_1+\eps}=y_1,\quad f_{\omega_2}=x_{12},\\
f_{\omega_3+\eps}=y_1 x_{23}-y_2 x_{13}+y_3 x_{12},\quad f_{\omega_4}=\pf x=x_{12} x_{34}-x_{13} x_{24}+x_{14} x_{23}
\end{gather*}
define four prime divisors $D_1$, $D_2$, $D_3$ and $D_4$ in $X$, respectively. The divisor $D_4$ is $G$-stable, and the others are the colors: $\mathcal{D}=\{D_1,D_2,D_3\}$ (the polynomial $f_{\omega_3+\eps}$ is the coefficient of $e^*_1\wedge e^*_2\wedge e^*_3$ in the expansion of the tensor $y\wedge x$ in the standard basis). Finally, it is easy to compute that $\mathcal{D}(\alpha_1)=\{D_1\}$, $\mathcal{D}(\alpha_2)=\{D_2\}$ and $\mathcal{D}(\alpha_3)=\{D_3\}$.

The variety $X$ has two closed $G$-orbits: $Y_1=\PP((\CC^4)^*)$, and the orbit $Y_2$ formed by the images in $\PP(V)$ of the decomposable bivectors (those that are exterior products of two vectors). For the colors whose closures contain $Y_1$ or $Y_2$ we have $\mathcal{F}_1=\{D_2,D_3\}$ and $\mathcal{F}_2=\{D_1,D_3\}$, and the cones $\mathcal{C}_1$ and $\mathcal{C}_2$ are generated by $v_{D_4}$ and the vectors $v_D$ with $D\in\mathcal{F}_i$. We describe the blowup of $Y_1$. Figure~\ref{fig: interesting_example}(a) shows the colored fan of the simple spherical subvariety $X_{Y_1}\subset X$ with the unique closed $G$-orbit $Y_1$ (we do not draw the color $D_1$, which does not occur in this fan, and show only the walls of the cones $\mathcal{C}_1$ and $\mathcal{V}$). In the basis dual to $\beta_1,\beta_2,\gamma$ we have $v_{D_2}=(1,1,-2)$, $v_{D_3}=(1,-1,1)$ and $v_{D_4}=(-1,0,0)$.

Let us now describe the colored fan $\FF(X')$ of the blowup. In the notation of Theorem~\ref{thm: main} we have $S_{\mathcal{F}_1}=\{\alpha_2,\alpha_3\}$, hence $\Sigma_*=\{\beta_1\}$, and the localized valuation cone $\mathcal{V}_*$ is the half-space $\{v\mid\pairing{v}{\beta_1}\leq0\}$. The intersection $\mathcal{C}_1\cap l\mathcal{V}_*$ is the cone generated by the vectors $v_1=(0,1,-2)=v_{D_2}+v_{D_4}$ and $v_2=(0,-1,1)=v_{D_3}+v_{D_4}$. Thus $\FF(X')=\Star(\FF(X),v_0)=\Star((\mathcal{C}_1,\mathcal{F}_1),v_0)\cup(\FF(X)\backslash\FF(X_{Y_1}))$, where $v_0=v_1+v_2=(0,0,-1)$. Explicitly,
\begin{align*}
\Star((\mathcal{C}_1,\mathcal{F}_1),v_0)=\{&(\cone(v_{D_3},v_{D_4},v_0),\{D_3\}),\\
&(\cone(v_{D_3},v_{D_4}),\{D_3\}),\\
&(\cone(v_{D_4},v_0),\varnothing),\\
&(\cone(v_{D_4}),\varnothing),\\
&(\cone(v_0),\varnothing),\\
&(0,\varnothing)\}.
\end{align*}
Figure~\ref{fig: interesting_example}(b) illustrates $\mathcal{V}_*$ and the cone $\mathcal{C}_1$.
\end{example}

\end{document}